\input amstex
\documentstyle{amsppt}
\catcode`\@=11
\def\logo@{}
\catcode`\@=\active \magnification\magstep1 \hsize 139mm
\addto\tenpoint{\normalbaselineskip12pt\normalbaselines}
\mathsurround1pt \tolerance500

\NoRunningHeads
\pageno 1
\TagsOnLeft


\define\N{\Bbb N}
\define\C{\Bbb C}

\define\a={\equiv}

\define\1{\bold 1}

\fontdimen14\tensy=3pt
\fontdimen16\tensy=3pt
\fontdimen17\tensy=3pt

\document

\topmatter
\title {Remarks on generalized Ramanujan sums and even functions}
\endtitle
\author {L\'aszl\'o T\'oth}
\endauthor
\affil {University of P\'ecs, Institute of Mathematics and Informatics,\\
Ifj\'us\'ag u. 6, H-7624 P\'ecs, Hungary \\
E-mail: {\tt ltoth\@ttk.pte.hu}\\ \ \\ {\it Acta Math. Acad.
Paedagog. Nyh\'azi. (N. S.), electronic, {\bf 20} (2004), 233-238}}
\endaffil
\abstract We prove a simple formula for the main value of $r$-even functions and give applications of it. Considering the generalized Ramanujan sums $c_A(n,r)$ involving regular systems $A$ of divisors we show that it is not possible to develop a Fourier theory with respect to $c_A(n,r)$, like in the the usual case of classical Ramanujan sums $c(n,r)$.
\endabstract
\endtopmatter



Supported partially by the Hungarian National Foundation for
Scientific Research under grant OTKA T031877.

\head{1. Introduction}
\endhead

Ramanujan's trigonometric sum $c(n,r)$ is defined as the sum of $n$-th powers of the $r$-th primitive roots of unity, i. e.
$$
c(n,r)=\sum \Sb k \text{(mod $r$)}\\
(k,r)=1\endSb \exp(2\pi i kn/r),
$$
where $r,n\in \N \equiv \{1,2,3,...\}$. In his original paper [R18] S. Ramanujan proved, among others, that for every $n\in \N$,
$$
\frac{\sigma(n)}{n}= \frac{\pi^2}{6} \sum_{r=1}^{\infty}
\frac{c(n,r)}{r^2}= \frac{\pi^2}{6} \left(1+\frac{(-1)^n}{2^2}+
\frac{2\cos(2\pi n/3)}{3^2} + \frac{2\cos(\pi n/2)}{4^2}+\cdots \right),
\tag1
$$
where $\sigma(n)$ stands for the sum of the positive divisors of $n$. Formula (1)
shows how the values of $\sigma(n)/n$ fluctuate harmonically about their mean value
$\pi^2/6$. Here the main value of a function $f:\N \to \C$ is defined by
$M(f)=\lim_{x\to \infty} \frac1{x} \sum_{n\le x} f(n)$, if this limit exists.

The orthogonality relations
$$
M(c(\cdot,r)c(\cdot,s))=\delta_{r,s}\phi(r), \tag2
$$
where $\delta_{r,s}$ is the Kronecker-symbol and $\phi(r)=c(r,r)$ is Euler's arithmetical function, suggest to have expansions, convergent pointwise
or in other sense, of functions $f$ of the form
$$
f(n)=\sum_{r=1}^{\infty} a_r c(n,r), \quad n\in \N,
$$
where the coefficients $a_r$ are given by
$$
a_r=\frac1{\phi(r)} M(fc_r).
$$

A Fourier analysis of arithmetical functions, with respect to Ramanujan sums, parallel to periodic and almost periodic functions, was developed
by several authors, cf. J. Delsarte [D45], W. Schwarz and J. Spilker [SchS74], A. Hildebrand [H84], G. G\'at [G91], L. Lucht [L95], see also the books [K75], [SchS94].

Ramanujan's sum $c(n,r)$ is an example of an $r$-even function (even function (mod $r$)), i.e. a function $f(n,r)$ such that $f(n,r)=f(gcd(n,r),r)$ for every $n,r\in \N$. This concept was introduced by E. Cohen [C55]. For a fixed $r$ the set $\Cal E_r$ of $r$-even functions $f:\N \to \C$ is a complex Hilbert space of finite dimension
$\tau(r)=$ the number of divisors of $r$, under the inner product
$$
\langle f,g \rangle = \frac1{r} \sum_{d|r} \phi(d)f(r/d)\overline{g(r/d)},
$$
and $(c'(\cdot,q))_{q|r}$, $c'(n,q)=\frac1{\sqrt{\phi(q)}}c(n,q)$ is an orthonormal basis for $\Cal E_r$. The main value $M(c(\cdot,r))$ exists for every $r\ge 1$, it is $M(c(\cdot,r))=\delta_{r,1}$, hence $M(f)$ exists for each $f\in \Cal E_r$.
If $\Cal E= \cup_{r\in \N} \Cal E_r$, then
$\Cal E$ is a dense subalgebra of the algebra $\C^{\N}$.

For these and various other properties of $c(n,r)$ and of $r$-even functions
see [K75], [McC86], [SchS94].

The Ramanujan sum $c(n,r)$ has been generalized in several directions. One of the generalizations, due to P. J. McCarthy [McC68], notation $c_A(n,r)$, is involving regular systems $A$ of divisors, see Section 2, and it has all
nice algebraic properties of the usual kind.

The following question can be formulated. Is it possible to develop a Fourier theory concerning the generalized sums $c_A(n,r)$, analogous to the usual one ?

The aims of this paper are the following:

- to prove a simple formula for the main value of $r$-even functions (Proposition 1), this result seems to have not been appeared in the literature, and to give applications of it,

- to compute the main value of $c_A(\cdot,r)$ for an arbitrary regular
system $A$ (Proposition 2),

- to show that the answer is negative for the question formulated above (Propositions 3 and 4).

\head{2. Regular convolutions}\endhead

Let $A(n)$ be a subset of the set of positive divisors of $n$ for
each $n\in \N$. The $A$-convolution
of the functions $f,g: \N \to \C$ is given by
$$
(f*_Ag)(n)=\sum_{d\in A(n)} f(d)g(n/d), \quad n\in \N.
$$

The system $A=(A(n))_{n\in \N}$ of divisors is called regular, cf. [N63], if

(a) $(\C^{\N},+,*_A)$ is a commutative ring with unity,

(b) the $A$-convolution of multiplicative functions is multiplicative
(recall that function $f$ is multiplicative if $f(mn)=f(m)f(n)$ whenever
$gcd(m,n)=1$),

(c) the constant $1$ function has
an inverse $\mu_A$ (generalized M\"obius function)
with respect to $*_A$ and $\mu_A(p^a)\in \{-1,0\}$ for every prime
power $p^a (a\ge 1)$.

It can be shown, cf. [N63], [McC86], that $*_A$ is regular iff

(i) $A(mn)=\{de: d\in A(m), e\in A(n)\}$ for every $m,n\in \N, (m,n)=1$,

(ii) for every prime power $p^a (a\ge 1)$ there exists a divisor
$t=t_A(p^a)$ of $a$, called the type of $p^a$ with respect
to $A$, such that $A(p^{it})=\{1,p^t,p^{2t},...,p^{it}\}$ for every $i\in
\{0,1,...,a/t\}$.

Examples of regular systems of divisors are $A=D$, where $D(n)$ is the set of all positive disors of $n$,  and $A=U$, where $U(n)$ is the set of divisors $d$ of $n$ such that $(d,n/d)=1$ (unitary divisors). For every prime power $p^a$ one has $t_D(p^a)=1$ and $t_U(p^a)=a$.
Here $*_D$ and $*_U$ are the Dirichlet convolution and the unitary convolution, respectively. For properties of regular convolutions and related arithmetical functions we refer to [N63], [McC86], [S78], [T97].

The following generalization of $c(n,r)$ is due to P. J. McCarthy
[McC68], see also [McC86]. For a regular system $A$ of divisors and $r,n\in \N$ let
$$
c_A(n,r)=\sum \Sb k \text{ (mod $r$) }\\
(k,r)_A=1\endSb \exp(2\pi i kn/r),
$$
where $(k,r)_A=\max \{d\in \N: d|k, d\in A(r)\}$, and let $c_A(r,r)\equiv \phi_A(r)$ be the generalized Euler function.
For $A=U$ the functions $c_U(n,r)\equiv c^*(n,r)$ and $\phi_U(r)\equiv \phi^*(r)$ were introduced by E. Cohen [C60].

$c_A(n,r)$ preserves the basic properties of $c(n,r)$. For example,
for every regular $A$ and $r, n\in \N$ one has
$$
c_A(n,r) =\sum_{d|n, d\in A(r)} d\mu_A(r/d), \tag3
$$
hence $c_A(n,r)$ is integer-valued and it is multiplicative in $r$. Note that
$$
c_A(n,r) =\sum_{d|r, \gamma_A(r)|d} c(n,d), \tag4
$$
where $\gamma_A$ is multiplicative and $\gamma_A(p^a)=p^{a-t+1}$ for every prime power $p^a$ ($a\ge 1$), here $t=t_A(p^a)$ (generalized core function), see [McC68], Th. 2.

The function $f$ is called $A$-even (mod $r$) if $f(n,r)=f((n,r)_A,r)$ for each $n,r\in \N$, cf. [McC68].
Let $\Cal E_{A,r}$ denote the set of functions $f(n,r)$ which are $A$-even (mod $r$). Then $\Cal E_{A,r} \subset \Cal E_r$ for every $A$ and $r\in \N$.
For example, $c_A(n,r)$ is $A$-even (mod $r$).
Let $\Cal E_A= \cup_{r\in \N} \Cal E_{A,r}$.

\head{3. Results}\endhead

\proclaim{Proposition 1} Let $r\in \N$ and $f\in \Cal E_r$. Then
$$
M(f)=\frac1{r} (f*_D\phi)(r) \equiv \frac1{r} \sum_{e|r} f(e)\phi(r/e).
$$
Moreover, for every $x\ge 1$ and every $\varepsilon >0$,
$$
\frac1{K_f} \left| \sum_{n\le x} f(n)- M(f)x \right| \le  C_{\varepsilon} r^{1+\varepsilon},
$$
where $|f(n)|\le K_f, n\in \N$ and $C_{\varepsilon}$ is a constant depending
only on $\varepsilon$.
\endproclaim

\demo{Proof} Write $f$ in the form
$$
f(n)=\sum_{q|r} h(q)c(n,q), \quad n\in \N,
$$
where the Fourier coefficients $h(q)$, are given for every $q|r$ by
$$
h(q)= \frac1{r\phi(q)} \sum_{e|r} \phi(e)f(r/e)c(r/e,q)=
\frac1{r}\sum_{e|r} f(r/e)c(r/q,e),
$$
and use that $\sum_{n\le x} c_q(n)=\delta_{q,1} x +R_q(x)$, where $|R_q(x)|\le q^{1+\varepsilon}$, cf. [McC86], Ch. 2, [K75], Ch. 7.
Then
$$
\sum_{n\le x} f(n)= \sum_{n\le x} \sum_{q|r} h(q)c(n,q) =
\sum_{q|r} h(q) \sum_{n\le x} c(n,q)=
$$ $$
=\sum_{q|r} h(q) (\delta_{q,1}x + R_q(x))= h(1)x + \sum_{q|r} h(q)R_q(x),
$$
where $h(1)=\frac1{r}(f*_D\phi)(r)$, $|h(q)|\le \frac1{r}\sum_{e|r} |f(r/e)||c(r/q,e)|\le K_f \frac1{r}\sum_{e|r} e =K_f\sigma(r)/r$ and
$$
|\sum_{q|r} h(q) R_q(x)| \le K_f (\sigma(r)/r) \sum_{q|r} q^{1+\varepsilon}
$$
and the result follows by the usual estimates $\sigma(r)\le r\tau(r)$ and
$\tau(r)\ll r^{\varepsilon}$.
\enddemo

As an application, consider the function $\phi(s,d,n)$ defined by $\phi(s,d,n)=\# \{k\in \N \cap [1,n]: (s+(k-1)d,n)=1\}$,
where $s,d\in \N, (s,d)=1$. Note that $\phi(1,1,n)=\phi(n)$ is the Euler function. T. Maxsein [M90] pointed out that $\phi(s,\cdot,n)$ is an $n$-even function and determined its main value:
$$
M(\phi(s,\cdot,n))=n\prod_{p|n} (1-1/p+1/p^2),
$$
where the product is over the prime divisors $p$ of $n$. This follows from Proposition 1 by easy computations.

Proposition 1 applies also for $f(n)=c_A(n,r)$, which is an $r$-even
function. However, a better error term can be obtained and the computations
are simpler by a direct proof using representation (3). We have

\proclaim{Proposition 2} For every regular system $A$ and $r\in \N$,
$$
M(c_A(\cdot,r))=\delta_{r,1}
$$
and
$$
|\sum_{n\le x} c_A(n,r)- \delta_{r,1}x|\le \psi_A(r),
$$
where $\psi_A$ is multiplicative and $\psi_A(p^a)=p^a+p^{a-t}$ for every prime power $p^a$ ($a\ge 1$), where $t=t_A(p^a)$ (generalized Dedekind function).
\endproclaim

\demo{Proof} Using (3),
$$
\sum_{n\le x} c_A(n,r)= \sum \Sb n\le x \\d|n, d\in A(r)\endSb d\mu_A(r/d)=
\sum_{d\in A(r)} d\mu_A(r/d) \sum \Sb n\le x\\ d|n \endSb 1=
$$
$$
=\sum_{d\in A(r)} d\mu_A(r/d) [x/d]=
x\sum_{d\in A(r)} \mu_A(r/d) - \sum_{d\in A(r)} d\mu_A(r/d)(x/d -[x/d])
=x\delta_{r,1}+R_A(r),
$$
where
$$
|R_A(r)|\le \sum_{d\in A(r)} d|\mu_A(r/d)|=\psi_A(r).
$$
\enddemo

Note that $\psi_A(r)\le \sigma(r)< Cr\ln \ln r$ for every $r\in \N$, with a suitable constant $C$.

The following result shows that for every system $A\ne D$ the orthogonality relations (2) are violated.

\proclaim{Proposition 3} For every regular system $A$,
$$
M(c_A(\cdot,r)c_A(\cdot,s))=\cases \phi_A(r), &\text{ if } r=s,\\
                             0, &\text{ if } rs>1, (r,s)=1,
\endcases
$$
but for $A\ne D$ there exist $r,s$ such that $r\ne s$ and $M(c_A(\cdot,r)c_A(\cdot,s))\ne 0$.
\endproclaim

\demo{Proof} Let $A$ be arbitrary. Applying (4) and (2) we obtain
$$
M(c_A(\cdot,r)c_A(\cdot,s))= \sum \Sb d|r, \gamma_A(r)|d \\ e|s, \gamma_A(s)|e \endSb M(c(\cdot,r)c(\cdot,s))=
\sum \Sb d|r, \gamma_A(r)|d \\ d|s, \gamma_A(s)|d \endSb \phi(d). \tag5
$$
Using that $\gamma_A(k)> 1$ for $k>1$ we get the first part of the desired result.

Now let $A\ne D$. Then there exists a prime power $p^a$ such that $t\equiv t_A(p^a)>1$. Therefore $A(p^t)=\{1,p^t\}$. Let $r=p$, $s=p^t$, then
$r\ne s$ and the last sum in (5) has one single term, namely $\phi(p)=p-1
\ne 0$.
\enddemo

\proclaim{Proposition 4} $\Cal E_A$ is a vector space if and only if $A=D$.
\endproclaim

\demo{Proof} We show that $\Cal E_A$ is not a vector space for $A\ne D$.

Suppose that $A\ne D$. Then there exists a prime power $p^a$ such
that $t\equiv t_A(p^a)>1$. Hence $A(p^t)=\{1,p^t\}$. Let
$$
f(n)=(n,p)_A= \cases &p, \text{ if } p|n,\\
                     &1, \text{ otherwise}, \endcases
$$ $$
g(n)=(n,p^t)_A= \cases &p^t, \text{ if } p^t|n, \\
                       &1, \text{ otherwise}. \endcases
$$
Then $f, g \in \Cal E_A$ and suppose that $h=f+g\in \Cal E_A$,
i. e. there exists $r\in \N$ such that $f+g\in \Cal E_{A,r}.$ Here
$$
h(n)=f(n)+g(n)= \cases &p+p^t, \text{ if } p^t|n, \\
                       &1+p, \text{ if } p|n, p^t\not| n, \\
                       &2, \text{ if } p\not| n. \endcases
$$
From $1+p=h(p)=h((p,r)_A)$ we have $(p,r)_A=p$, $p\in A(r)$ and
from $p+p^t=h(p^t)=h((p^t,r)_A)$ we obtain $(p^t,r)_A=p^t$, $p^t\in A(r)$.

Let $r=p^k s$, where $p\not| s$. Then $p\in A(p^k)$ and $p^t\in A(p^k)$,
therefore $p\in A(p^t)$, in contradiction with $A(p^t)=\{1,p^t\}$.
\enddemo

\Refs

[C55] {\smc E. Cohen}, A class of arithmetical functions, Proc. Nat. Acad. Sci. U.S.A. {\bf 41} (1955), 939-944.

[C60] {\smc E.Cohen}, Arithmetical functions associated with the unitary
divisors of an integer, Math. Z. {\bf 74} (1960), 66-80.

[D45] {\smc J. Delsarte}, Essai sur l'application de la th\'eorie des fonctions presque-p\'eriodiques \'a l'arithm\'eti-\newline que, Ann. Sci. \'Ecole Norm. Sup. {\bf 62} (1945), 185-204.

[G91] {\smc G. G\'at}, On almost even arithmetical functions via orthonormal systems on Vilenkin groups, Acta Arith. {\bf 60} (1991), 105-123.

[H84] {\smc A. Hildebrand}, \"Uber die punktweise Konvergenz von Ramanujan-Entwiklungen zahlentheoretischer Funktionen, Acta Arith. {\bf 44} (1984), 109-140.

[K75] {\smc J. Knopfmacher}, Abstract Analytic Number Theory, North Holland
Publ. Co., Amsterdam - Oxford, 1975.

[L95] {\smc L. Lucht}, Weighted relationship theorems and Ramanujan expansions, Acta Arith. {\bf 70} (1995), 25-42.

[M90] {\smc T. Maxsein}, A note on a generalization of Euler's
$\varphi$-function, Indian J. Pure Appl. Math. {\bf 21} (1990), 691-694.

[McC68] {\smc P. J. McCarthy}, Regular arithmetical convolutions,
Portugal. Math. {\bf 27} (1968), 1-13.

[McC86] {\smc P. J. McCarthy}, Introduction to Arithmetical Functions,
Springer- Verlag, New York - Berlin - Heidelberg - Tokyo, 1986.

[N63] {\smc W. Narkiewicz}, On a class of arithmetical convolutions,
Colloq. Math. {\bf 10} (1963), 81-94.

[R18] {\smc S. Ramanujan}, On certain trigonometric sums and their applications in the theory of numbers, Transactions Cambridge Phil. Soc.
{\bf 22} (1918), 259-276.

[SchS74] {\smc W. Schwarz, J. Spilker}, Mean values and Ramanujan expansions of almost even arithmetical functions, in Topics in Number Theory, Colloq. Math. Soc. J. Bolyai {\bf 13} (1974), 315-357.

[SchS94] {\smc W. Schwarz, J. Spilker}, Arithmetical Functions,
London Math. Soc. Lecture Notes Series {\bf 184}, Cambridge Univ. Press,
1994.

[S78] {\smc V. Sita Ramaiah}, Arithmetical sums in regular
convolutions, J. Reine Angew. Math. {\bf 303/304} (1978), 265-283.

[T97] {\smc L. T\'oth}, Asymptotic formulae concerning arithmetical
functions defined by cross-convolutions, I. Divisor-sum functions and
Euler-type functions, Publ. Math. Debrecen {\bf 50} (1997), 159-176.

\endRefs

\enddocument